\documentclass{article}
\usepackage{amsmath, enumerate}
\title{Mystery of the $18$th Elephant Solved}
\author{{\bf Dr V N Krishnachandran}\\ Vidya Academy of Science \& Technology\\ Thrissur - 680501, Kerala, India \\
email: {\tt  krishnachandran@vidyaacademy.ac.in}}
\date{}
\begin{document}
\maketitle

The {\em The 18th Elephant -- Three Monologues} is a chilling 63-minute documentary film  depicting  the insufferable cruelty handed out to the biggest mammals both in the wild and in the domesticated habitats.
 The film was produced by Alternate Network of Media People (ANMPE), a non-profit Indian collective for exploring the possibilities of the media for advancing the cause of social justice, environmental protection and conservation, and it went on to receive several national and international awards. The producer and the director were awarded the Best Scientific Film/Best Environment/Conservation/Preservation Film by Govt of India as part of the $51$st National Film Awards announced in February 2005.
{\em The $18$th Elephant} had also been selected for the prestigious Green Oscar Festival held in Bristol in October 2004 and secured the Green Oscar-Panda Award.
The Green Oscar, considered the biggest environmental film festival, is organised by the London-based Ashden Trust in collaboration with the British government.

But our focus here is not the plight of the elephants, but the mystery behind the $18$th elephant!

In a well known legend, a wise man introduces an elephant, which happens to be the $18$th elephant, mysteriously to help solve a mathematical conundrum and then disappears once a solution was obtained to the problem at hand.  
\section{The story of the $18$th elephant}
The story goes something like this (see, for example, \cite{IdeaTransfer}):
\begin{quote}
This is a story about a wealthy man in ancient India who owned a herd of elephants. As times passed by, he became old he decided to pass on his wealth to his three sons.  
Because his sons had very different natures and different experience, he decided to divide the assets accordingly.  And, to prevent any disputes from arising he determined to make the terms of the legacy irrevocable. He declared that he would divide his elephants based on their skills and attitudes. 
One half of his elephants would go to his eldest son, one third to his middle son, and one-ninth to his youngest son.  At first everyone seemed to agree that this was both prudent and fair.  But then a buzz began and turned into arguments. The company owned seventeen elephants.  What was half of seventeen?  What was one third?  One ninth?  

What could the wealthy man do?  His decision was irrevocable.  He was bewildered and demoralized.  He fled to his guru, fell on his knees, and begged for a solution.  The guru listened and smiled.  There was no problem, he told the sishya.  The guru himself had an elephant.  Just take the elephant back to town and add it to his own to complete the division.  After that, return the elephant to the guru.

He was speechless.  Was his guru that simple?  What could this additional elephant accomplish?  But he did as the guru advised, and everyone gathered round the elephants to watch.

With eighteen elephants the oldest son claimed nine elephants for him, leaving nine for his brothers to divide.  The second son claimed his third and took six elephants, leaving three.  The youngest son took his one-ninth portion, namely, two elephants. All the three sons among themselves took a total of nine plus six plus  two, that is seventeen, elephants only. 

This left only the guru’s elephant, which the man happily returned to the guru.
\end{quote}

The exact source of the story could not be located. According to the producers of {\em The $18$th Elephant}, it was a story told by Budha (see \cite{Hindu}). The story has perhaps universal appeal, for its essence  can be seen rendered  in different ways. For example, there is a story involving seventeen horses to be divided among three sons exactly as in the story of the the $18$th elephant (see \cite{Change}).

The story is generally presented to teach the moral: ``The attitude of negotiation and problem solving is to find the 18th elephant, that is,  the common ground. Once a person is able to find the 18th elephant the issue is resolved. It is difficult at times. However, to reach a solution, the first step is to believe that there is a solution. If we think that there is no solution, we will not be able to reach any!'' (See, for example, \cite{Nuggets}.)
\section{The arithmetic of the $18$th elephant}
Let us now return to the mathematical mystery of the $18$th elephant. Recall that the total number of elephants in the herd was $17$. 
\begin{itemize}
\item
The first son is eligible to receive one-half of the total number of elephants, that is, $\tfrac{17}{2}=8\frac{1}{2}$ elephants.
\item 
The second son is eligible to receive one-third of the total number of elephants, that is, $\tfrac{17}{3}=5\frac{2}{3}$ elephants.
\item
The third son is eligible to receive one-ninth of the total number of elephants, that is, $\tfrac{17}{9}=1\frac{8}{9}$ elephants.
\end{itemize}
The total number of elephants the sons are to receive is 
$$
8\frac{1}{2} + 5\frac{2}{3} + 1\frac{8}{9} = 16\frac{1}{18}.
$$
This sum is not equal to $17$, and is in fact less than $17$. What this implies is that even if the sons are given elephants as wished by the father some of the elephants, namely,
$$
17 - 16\frac{1}{18} = \frac{17}{18}
$$
will be left out. The wise guru exploited this to solve the problem. The guru decided to divide the balance of $\frac{17}{18}$ elephants according to his fancy to help the father solve the problem. The guru, without telling so explicitly to the father, was deciding to give the $\frac{17}{18}$ elephants to the three sons as follows:
\begin{itemize}
\item
Give $\frac{1}{2}$ elephant to the first son so that he would get a total of $8\frac{1}{2}+\frac{1}{2}=9$ elephants.
\item
Give $\frac{1}{3}$ elephant to the second son so that he would get a total of $5\frac{2}{3}+\frac{1}{3}=6$ elephants.
\item
Give $\frac{1}{9}$ elephant to the third son so that he would get a total of $1\frac{8}{9}+\frac{1}{9}=2$ elephants.
\end{itemize}
Fortuitously for the guru, the additional elephants given to the sons add up to the balance of elephants remaining after the herd had been divided as per the wish of the father:
$$
\frac{1}{2}+\frac{1}{3} + \frac{1}{9}=\frac{17}{18}.
$$
Also the numbers of the elephants obtained by the three sons turned out to be exactly one-half, one-third and one-ninth of $18$!

All this worked out the way it did because of certain arithmetical relations among the numbers $2$, $3$, $9$, $17$ and $18$ as we shall see shortly.
\section{The mathematics of the $18$th elephant}
\subsection{A similar situation}
A man has $57$ elephants and has four sons. He wishes to give his sons, respectively, one-third, one-fourth, one-fifth and one-sixth of the total number of elephants. Obviously, this cannot be done. This time, the guru will have to loan three elephants for the computations. With the three elephants obtained on loan from the guru, the total comes to $60$. The sons will get respectively $20$, $15$, $12$ and $10$ elephants. These add up to $57$. The elephants obtained on loan from the guru can be returned as unused.  
\subsection{Now the mathematics}
Let there be $N$ elephants in the herd and let it be required to divide the herd into groups in the ratios 
$$
\frac{1}{s_1} : \frac{1}{s_2}: \cdots : \frac{1}{s_k}
$$
among the $k$ sons of the owner of the herd. It will be assumed that $s_1, s_2, \cdots, s_k$ are all integers. The shares the sons are to receive are $\frac{N}{s_1}$, $\frac{N}{s_2}$, $\cdots$, $\frac{N}{s_k}$ and these may not be all integers. Let the number of elephants received on loan from the guru be $x$ so that the total number of elephants now becomes $N+x$. For the wise guru's advise to work we must have
\begin{equation}\label{eq1}
\frac{N+x}{s_1}+\frac{N+x}{s_2}+\cdots + \frac{N+x}{s_k}=N.
\end{equation}
The following conclusions can be derived from this relation:
\begin{enumerate}
\item
We must have
$$
\frac{1}{s_1} + \frac{1}{s_2} + \cdots + \frac{1}{s_k}<1.
$$
\item
Obviously, the values of the fractions $\frac{N+x}{s_1}$, $\frac{N+x}{s_2}$, $\cdots$, $\frac{N+x}{s_k}$ must all be exact integers. This implies that $N+x$ must be a multiple of each of the $k$ integers $s_1, s_2, \cdots, s_k$. Putting this in another way, the least common multiple of $s_1, s_2, \cdots, s_k$ must be a divisor of $N+x$.
\item
Let us rewrite Eq.\eqref{eq1} in the following form:
\begin{equation}\label{eq2}
\frac{1}{s_1} + \frac{1}{s_2}+\cdots + \frac{1}{s_k}=\frac{N}{N+x}.
\end{equation}
Writing 
$$
m = \text{lcm } \{s_1, s_2, \ldots, s_k\}
$$
we can express the left hand side of Eq.\eqref{eq2} in the following form
$$
\frac{1}{s_1}+\frac{1}{s_2}+\cdots+\frac{1}{s_k}=\frac{r}{m}
$$
for some integer $r$.
\item
Since $N+x$ is a multiple of $m$, we can find an integer $a$ such that
$$
N+x=am.
$$
Hence, from Eq.\eqref{eq2}, we have
$$
\frac{r}{m}=\frac{N}{am}
$$
from which it follows that
$$
N=ar.
$$
\end{enumerate}
\subsubsection*{Conclusion}
In the notations introduced above, we may conclude that the guru's advice will work only if the number of the elephants in the herd, that is $N$, is a multiple of $r$. Also if $N=ar$, then the number of elephants to be obtained on from the guru is $am-N$.
\section{Illustration}
Assume that a herd of elephants  is to be divided into four groups in the following ratio:
$$
\frac{1}{3}:\frac{1}{6}:\frac{1}{9}:\frac{1}{12}.
$$
We have 
$$
\text{lcm }\{3,6,9,12\}=36
$$
and
$$
\frac{1}{3}+\frac{1}{6}+\frac{1}{9}+\frac{1}{12}=\frac{25}{36}.
$$
In the notations introduced above we have $r=25$ and $m=36$. The total number of elephants $N$ can be a multiple of $25$. We may take $N=50$ so that $a=2$. In this case, the guru will be required to supply 
$$
x=am-N=2\times 36 - 50 =22
$$
elephants! With the the guru's elephants added to the herd, the total comes out to be $72$. The sons will receive respectively $\frac{72}{3}=24$, $\frac{72}{6}=12$, $\frac{72}{9}=8$ and $\frac{72}{12}=6$ numbers of elephants the total of which is 
$50$. The 22 elephants taken on loan from the guru can now be safely returned.
\section{Mystery solved!}
And that solves the mystery of the $18$th elephant!

\end{document}